\documentclass[11pt]{article}
\usepackage[top=2.54cm, bottom=2.54cm, left=2.70cm, right=2.70cm]{geometry}
\usepackage[tbtags]{amsmath}
\usepackage{amssymb}
\usepackage{amsthm}
\usepackage{fancyhdr}
\usepackage{latexsym}
\usepackage{mathrsfs}
\usepackage{wasysym}
\usepackage{fancyhdr}
\usepackage{float}
\usepackage{graphicx}
\usepackage{tikz}
\usepackage[numbers,sort&compress]{natbib}
\usepackage{setspace}
\allowdisplaybreaks[4]

\newtheorem{theorem}{\bf Theorem}[section]
\newtheorem{lemma}[theorem]{Lemma}
\newtheorem{proposition}[theorem]{Proposition}

\newtheorem{problem}[theorem]{Problem}
\newtheorem{corollary}[theorem]{Corollary}
\newtheorem{conj}[theorem]{Conjecture}

\newtheorem{claim}{\indent Claim}[]

\begin{document}
\begin{spacing}{1.1}
\title{A new connectivity bound for a tournament to be highly linked}
\author{Bin Chen$^a$, \, Xinmin Hou$^{a,b,c}$, \, Gexin Yu$^{d}$, \, Xinyu Zhou$^b$\\
\small$^a$Hefei National Laboratory, Hefei 230088, Anhui, China\\
\small $^b$School of Mathematical Sciences,\\
\small  University of Science and Technology of China, Hefei 230026, Anhui, China\\
\small$^c$ CAS Key Laboratory of Wu Wen-Tsun Mathematics,\\
\small  University of Science and Technology of China, Hefei 230026, Anhui, China\\
\small $^d$Department of Mathematics, William $\&$ Mary, Williamsburg, VA, USA}
\date{}
\maketitle
{\bf Abstract}\,: A digraph $D$ is $k$-linked if for any pair of two disjoint sets $\{x_{1},x_{2},\ldots,x_{k}\}$ and $\{y_{1},y_{2},\ldots,y_{k}\}$ of vertices in $D$, there exist vertex disjoint dipaths $P_{1},P_{2},\ldots,P_{k}$ such that $P_{i}$ is a dipath from $x_{i}$ to $y_{i}$ for each $i\in[k]$. Pokrovskiy (JCTB, 2015) confirmed a conjecture of K\"{u}hn et al. (Proc. Lond. Math. Soc., 2014) by verifying that every $452k$-connected tournament is $k$-linked. Meng et al. (Eur. J. Comb., 2021) improved this upper bound by showing that any $(40k-31)$-connected tournament is $k$-linked. In this paper, we show a better upper bound by  proving that every $\lceil 12.5k-6\rceil$-connected tournament with minimum out-degree at least $21k-14$ is $k$-linked.
Furthermore, we improve a key lemma that was first introduced by Pokrovskiy (JCTB, 2015) and later enhanced by Meng et al. (Eur. J. Comb., 2021).

{\bf AMS}\,: 05C20; 05C38; 05C40.


\section{Introduction}
Let $D=(V,A)$ be a digraph with vertex set $V(D)$ and arc set $A(D)$. All digraphs we considered throughout this paper are finite and simple.
For a vertex $v\in V(D)$, the \emph{out-neighborhood} of $v$ is the set $N^{+}_{D}(v)=\{u\in V(D): v\rightarrow u\}$, and the \emph{out-degree} of $v$ is denoted by $d^{+}_{D}(v)=|N^{+}_{D}(v)|$.
Analogously, the \emph{in-neighborhood} of $v$ is the set $N^{-}_{D}(v)=\{u\in V(D) : u\rightarrow v\}$, and the \emph{in-degree} of $v$ is denoted by $d^{-}_{D}(v)=|N^{-}_{D}(v)|$.
By $\delta^{+}(D)$ and $\delta^{-}(D)$ we mean the \emph{minimum out-degree} and \emph{minimum in-degree} of $D$, respectively.
The {\em minimum semi-degree} $\delta_0(D)=\min\{\delta^+(D),\delta^-(D)\}$.

A \emph{directed path} or {\em dipath} for short  of length $t$ of $D$ is a list of $t+1$ distinct vertices $v_{0},v_{1},\ldots,v_{t}$ such that $(v_{i}, v_{i+1})\in A(D)$ for all $0\leq i\leq t-1$. Let $P_{uv}$ be a dipath from $u$ to $v$ and we denote $\ell(P_{uv})$ as its length. The vertices $u$ and $v$ linked by $P_{uv}$ are called its \emph{initial} and \emph{end}, respectively. Two dipaths $P$ and $Q$ are said to be \emph{vertex disjoint}\;(\emph{disjoint} for short) if $V(P)\cap V(Q)=\emptyset$.

A digraph $D$ is \emph{strongly connected}\;(\emph{connected} for short) if any two vertices $u,v\in V(D)$ are connected by a dipath $P_{uv}$. A digraph $D$ is said to be \emph{$k$-connected} if it remains connected after the removal of any set of at most $k-1$ vertices of $D$. A digraph is {\em $k$-linked} if for any two disjoint vertex sets $\{x_{1},x_{2},\ldots,x_{k}\}$ and $\{y_{1},y_{2},\ldots,y_{k}\}$, there are disjoint dipaths $P_{1},P_{2},\ldots,P_{k}$ such that $P_{i}$ is a dipath from $x_{i}$ to $y_{i}$ for every $i\in \{1,2,\ldots, k\}$.
The relations between connectedness and linkedness are somewhat different in digraphs.
For example,  there is a function $f(k)$ such that every $f(k)$-connected graph is $k$-linked (such a function $f(k)$ was first given by Larman and Mani \cite{LM} and Jung \cite{J}, and some improvements of $f(k)$ were made by  Bollob\'{a}s and Thomason~\cite{BT} and by Thomas and Wollan \cite{TW05}). However,  there is no function $f(k)$ such that every strongly $f(k)$-connected directed graph is $k$-linked.
Indeed, Thomassen~\cite{Tho91} constructed digraphs of arbitrarily high connectedness that are not even $2$-linked. This indicates that there is no function for a general digraph of high connectedness to have large linkedness. Consequently, scholars have investigated this problem in special classes of digraphs, such as tournaments.

A \emph{tournament} is a digraph in which there is exactly one arc between every pair of distinct vertices. Thomassen~\cite{T84} verified that there is a function $g(k)$ such that every $g(k)$-connected tournament is $k$-linked, where $g(k)\leq Ck!$ for some constant $C$.
It will be very interesting to determine the exact value of $g(k)$ for any given $k$.
\begin{problem}\label{PROB: g(k)}
Determining the exact value of $g(k)$ for each given $k$.
\end{problem}
Except for the trivial exact value $g(1)=1$, there are few known exact values of $g(k)$. Here we list some related results to the best of our knowledge.
\begin{itemize}
\item[(1)] Thomassen~\cite{T84} determined that $g(2)=5$,  and Bang-Jensen~\cite{JBJ} extended this result to semicomplete digraphs, where a  semicomplete digraph is a digraph with no nonadjacent vertices.

\item[(2)] K\"{u}hn, Lapinskas, Osthus, Patel~\cite{KLOP} greatly improved the upper bound $g(k)\leq Ck!$ by showing  that $g(k)\leq 10^{4}k\log k$. They further conjectured that $g(k)$ is linear in $k$;

\item[(3)] Pokrovskiy~\cite{Pok} confirmed the above conjecture by showing that $g(k)\leq 452k$.


\item[(4)]  Meng, Rolek, Wang, Yu~\cite{MRWY} reduced the upper bound of $g(k)$ to $40k-31$.

\item[(5)] Authors in \cite{BJJ22,ZhouQiYan23} made attempts to improve the above bounds, but there are gaps in their proofs. See the last section for detail explanation.
\end{itemize}
Pokrovskiy~\cite{Pok} also showed that low connectivity cannot guarantee $k$-linked by constructing an infinite family of $(2k-2)$-connected tournaments with at least $6k$ vertices and arbitrarily large semi-degree but are not $k$-linked. Moreover, Pokrovskiy conjectured the following:
\begin{conj}[Pokrovskiy~\cite{Pok}]\label{CONJ: Pok}
For every $k$, there is $d(k)$ such that any $2k$-connected tournament with semi-degree at least $d(k)$ is $k$-linked.
\end{conj}

There are few improvements to Conjecture~\ref{CONJ: Pok}. However, there are a few results on this flavor. The following is a list of them.
\begin{itemize}
\item[(5)] (Gir\~{a}o and Snyder~\cite{GS}) There exists $d(k)$ such that every $4k$-connected tournament with minimum out-degree at least $d(k)$ is $k$-linked.

\item[(6)] (Gir\~{a}o, Popielarz, and Snyder~\cite{GPS}) Each $(2k+1)$-connected tournament with minimum out-degree at least $Ck^{31}$ is $k$-linked for some constant $C$.
\end{itemize}
Gir\~{a}o, et al.,~\cite{GPS} further constructed an infinite family of $(2.5k-1)$-connected tournaments that are not $k$-linked, yielding that the condition on large minimum out-degree is necessary.
Inspired by these works, we prove a result of the flavor of Conjecture~\ref{CONJ: Pok} in this article.

\begin{theorem}\label{main theorem}
For any positive integer $k$, every $\lceil 12.5k-6\rceil$-connected tournament with minimum out-degree at least $21k-14$ is $k$-linked.
\end{theorem}

As a corollary, we have  a better upper bound of $g(k)$ than (4) given by Meng et al.~\cite{MRWY}.
\begin{corollary}
	
$$g(k)\le 21k-14.$$
	
\end{corollary}	

An important component in our proof is the following  `anchor' lemma, initially discovered in~\cite{Pok}, and later  improved in ~\cite{MRWY}.
We make a further improvement of this lemma, which may have its own interest.
Let $X=\{x_{1},x_{2},\ldots,x_{k}\}$ and $Y=\{y_{1},y_{2},\ldots,y_{k}\}$ be two disjoint vertex sets of a digraph $D$. If there exist $k$ disjoint dipaths $P_{x_{1}y_{\pi(1)}}, P_{x_{2}y_{\pi(2)}},\ldots, P_{x_{k}y_{\pi(k)}}$ for any permutation $\pi$ of $\{1,2,\ldots,k\}$, then we call $X$ \emph{anchors} $Y$ in $D$. 

\begin{lemma}[Anchor lemma]\label{key lemma}
	Let $k\geq 1$ be an integer and let $T$ be a tournament on $n$ vertices. If $n\geq 8.5k-6$, then there are two disjoint vertex sets each of size $k$ such that one anchors another.
\end{lemma}

The rest of this paper is organized as follows. In the next section, we introduce some notation and present several auxiliary results, including the proof of Lemma~\ref{key lemma}.
The proof of Theorem~\ref{main theorem} will be given in Section 3. We give some remarks and discussions in the last section.

\section{Preliminaries and Proof of Lemma~\ref{key lemma}}

\noindent


\subsection{Notation and Preliminaries}

\noindent

An arc of $D$ with \emph{tail} $u$ and \emph{head} $v$ is denoted by $(u,v)$, and we say that $u$ dominates $v$ or $v$ is dominated by $u$. In addition, we write  $u\rightarrow v$ for $(u,v)\in A(D)$ and $u\nrightarrow v$ for $(u,v)\notin A(D)$.
For integers $a$ and $b$, let $[a,b]=\{a, a+1, \dots, b\}$ and write $[k]$ for $[1,k]$ for convenience.

Let $S\subseteq V(D)$. Write $N^{+}_{D}(v, S)=\{u\in S : (v, u)\in A(D) \}$ and $d^{+}_{D}(v, S)=|N^{+}_{D}(v, S)|$.
We can similarly define $N^{-}_{D}(v, S)=\{u\in S : (u, v)\in A(D) \}$ and $d^{-}_{D}(v, S)=|N^{-}_{D}(v, S)|$.

For any $W\subseteq V(D)$, let $D[W]$ be the subdigraph induced by $W$ in $D$. We write $D-W$ for the digraph $D[V(D)\backslash W]$. Moreover, we denote $N^{+}_{D}(W)=\{u\in V(D)\backslash W: \exists\;w\in W\;\text{such}\;\text{that}\;w\rightarrow u\}$, $N^{-}_{D}(W)=\{u\in V(D)\backslash W: \exists\;w\in W\;\text{such}\;\text{that}\;u\rightarrow w\}$, respectively. For any disjoint $X,Y\subseteq V(D)$, denote $A_D(X,Y)=\{(x,y)\in A(D) : x\in X\;\text{and}\;y\in Y\}$.

Given an order $v_{1},v_{2},\ldots,v_{n}$ of a digraph $D$, we say that an arc $(v_{i},v_{j})\in A(D)$ is a \emph{forward arc} if $i<j$ and a \emph{backward arc} if $i>j$, respectively. A \emph{median order} of a digraph is a vertex order that maximizes the number of its forward arcs.
The following are two basic properties of median orders of tournaments\;(see e.g., \cite{BM} page 101 for details).
\begin{proposition}[\cite{BM}]\label{median order}
Let $T$ be a tournament and $(v_{1},v_{2},\ldots,v_{n})$ a median order of $T$. Then, for any $i,j$ with $1\leq i<j\leq n$, the following hold.

(i). The interval $(v_{i},v_{i+1},\ldots,v_{j})$ is a median order of $T[\{v_{i},v_{i+1},\ldots,v_{j}\}]$;

(ii). The vertex $v_{i}$ dominates at least half of the vertices in $\{v_{i+1},v_{i+2},\ldots,v_{j}\}$ and the vertex $v_{j}$ is dominated by at least half of the vertices in $\{v_{i},v_{i+1},\ldots,v_{j-1}\}$.

\end{proposition}



\begin{proposition} \label{fact}
Every tournament on $n$ vertices has minimum out-degree at most $(n-1)/2$.
\end{proposition}

The celebrated Menger's Theorem\;(directed version) \cite{M} tells us that a digraph is $k$-connected if and only if for any two different vertices $x,y$, there exist $k$ internally disjoint dipaths from $x$ to $y$. We shall use an immediate corollary of Menger's Theorem as shown below:

\begin{corollary} \textnormal{(\cite{M})} \label{corollary}
Let $D$ be a $k$-connected digraph. Then for any two disjoint sets $\{x_{1},,\ldots,x_{k}\}$ and $\{y_{1},,\ldots,y_{k}\}$ of vertices of $D$, there are disjoint dipaths $P_{x_{1}y_{\pi(1)}},\ldots,P_{x_{k}y_{\pi(k)}}$ for some permutation $\pi$ of $\{1,2,\ldots, k\}$.
\end{corollary}

We also need the following lemma.

\begin{lemma}[\cite{O}] \label{matching}
	Let $G$ be a bipartite graph with partition $V(G)=X\cup Y$. Then the maximum number of vertices of $X$ that can be covered by a matching of $G$ is $|X|-\text{max}(\{|S|-|N_{G}(S)|:S\subseteq X\})$.
\end{lemma}


\subsection{Proof of Lemma~\ref{key lemma}}

Let $V(T)=\{v_{1},v_{2},\ldots,v_{n}\}$, and, without loss of generality, assume that $(v_{1},v_{2},\ldots,v_{n})$ forms a median order of $T$. Let $X=\{x_{1},x_{2},\ldots,x_{k}\}$ and $Y=\{y_{1},y_{2},\ldots,y_{k}\}$ such that $x_{i}=v_{i}$ and $y_{i}=v_{n-i+1}$ for $i\in [k]$. We also denote $Z=\{z_{1},z_{2},\ldots,z_{n-2k}\}$ such that $z_{i}=v_{k+i}$ for any $i\in [n-2k]$. Note first that $X\cup Y\cup Z=V(T)$ and they are mutually disjoint since $n\geq 8.5k-6$. Moreover, we write $X_{i}$\;(resp., $Y_{i}$) for $\{x_{i},x_{i+1},\ldots,x_{k}\}$\;(resp., $\{y_{i},y_{i+1},\ldots,y_{k}\}$). By Proposition~\ref{median order}, the interval $(x_{i},x_{i+1},\ldots,x_{k},z_{1},z_{2},\ldots,z_{n-2k})$ is a median order of the induced tournament $T[X_{i}\cup Z]$ and $d^{+}_{T[X_{i}\cup Z]}(x_{i})\geq (k-i+n-2k)/2=(n-k-i)/2$. This yields that$$d^{+}_{T[Z]}(x_{i})=d^{+}_{T[X_{i}\cup Z]}(x_{i})-d^{+}_{T[X_{i}]}(x_{i})\geq (n-k-i)/2-(k-i)=(n-3k+i)/2.$$
By symmetry, we can also deduce that$$d^{-}_{T[Z]}(y_{i})=d^{-}_{T[Y_{i}\cup Z]}(y_{i})-d^{-}_{T[Y_{i}]}(y_{i})\geq (n-k-i)/2-(k-i)=(n-3k+i)/2.$$

If $X$ anchors $Y$, then we are done. Therefore, we may assume that $X$ does not anchor $Y$, implying that there is some permutation $\pi$ of $[k]$ such that some desired dipaths from $x_{i}$ to $y_{\pi(i)}$ do not exist. Let $h$ be the maximal number such that $P_{1},P_{2},\ldots,P_{h}$ are disjoint dipaths in $T$ satisfying that for any $j\in [h]$:
\begin{itemize}
\item[(i)] if $j\equiv 1$\;(\text{mod} 2), then $P_{j}$ is a dipath $P_{x_{i}y_{\pi(i)}}$ such that $x_{i}\in X\backslash (V(P_{1})\cup V(P_{2})\cup \ldots \cup V(P_{j-1}))$ for minimal index $i$;

\item[(ii)] if $j\equiv 0$\;(\text{mod} 2), then $P_{j}$ is a dipath $P_{x_{\pi^{-1}(i)}y_{i}}$ such that $y_{i}\in Y\backslash (V(P_{1})\cup V(P_{2})\cup \ldots \cup V(P_{j-1}))$ for minimal index $i$;

\item[(iii)] every $P_{j}$ has length at most 3;

\item[(iv)] every $P_{j}$ has no interior vertex in $X\cup Y$.
\end{itemize}

Obviously, $h\leq k-1$. For simplicity, write $V(\widetilde{P}_{h})$ for $V(P_{1})\cup V(P_{2})\cup \ldots \cup V(P_{h})$. Here to proceed with the proof, we will divide the discussions into the following two cases, i.e., $h\equiv 0$\;(\text{mod} 2) and $h\equiv 1$\;(\text{mod} 2).

\vspace{0.2cm}

\noindent{\bf Case 1}. $h\equiv 0$\;(\text{mod} 2).
\vspace{0.2cm}

Let $i$ be the minimal index such that $x_{i}\in X\backslash V(\widetilde{P}_{h})$. By (i) and (ii), we have $i\geq h/2+1$ and $\pi(i)\geq h/2+1$.
As every $P_{j}$ has length at most three for each $j\in [h]$, we obtain that $x_{i}$\;(resp., $y_{\pi(i)}$) has at most two out-neighbors\;(resp., in-neighbors) in $V(P_{j})\backslash (X\cup Y)$. We write $X^{*}$ and $Y^{*}$ for $N^{+}_{T[Z\backslash V(\widetilde{P}_{h})]}(x_{i})$ and $N^{-}_{T[Z\backslash V(\widetilde{P}_{h})]}(y_{\pi(i)})$, respectively.
Then, we have
$$|X^{*}|\geq (n-3k+i)/2-2h\geq (n-6.5k+4.5)/2\geq k-1,$$and$$|Y^{*}|\geq (n-3k+\pi(i))/2-2h\geq (n-6.5k+4.5)/2\geq k-1.$$
By the maximality of $h$, there is no dipath of length at most three from $x_i$ to $y_{\pi(i)}$. Therefore, we see that $X^{*}\cap Y^{*}=\emptyset$ and $A(X^{*}\cup \{x_{i}\},Y^{*}\cup \{y_{\pi(i)}\})=\emptyset$.
Note that $|X^{*}\cup \{x_{i}\}|\geq k$ and $|Y^{*}\cup \{y_{\pi(i)}\}|\geq k$.
We can choose $Y^{**}\subseteq Y^{*}\cup \{y_{\pi(i)}\}$ and $X^{**}\subseteq X^{*}\cup \{x_{i}\}$ such that $|Y^{**}|=|X^{**}|=k$ and $y\rightarrow x$ for every $y\in Y^{**}$ and $x\in X^{**}$, i.e., $Y^{**}$ anchors $X^{**}$, as desired.

\vspace{0.2cm}

\noindent{\bf Case 2}. $h\equiv 1$\;(\text{mod} 2).
\vspace{0.2cm}

Let $i$ be the minimal index with $y_{i}\in Y\backslash V(\widetilde{P}_{h})$. Again by (i) and (ii),  we have $i\geq (h-1)/2+1=(h+1)/2$ and $\pi^{-1}(i)\geq (h+1)/2+1=(h+3)/2$. Since every $P_j$ has length at most three,  we have that $x_{\pi^{-1}(i)}$\;(resp., $y_{i}$) has at most two out-neighbors \;(resp., in-neighbors) in $V(P_{j})\backslash (X\cup Y)$ for $j\in [h]$. Write $X^{'}$ and $Y^{'}$ for $N^{+}_{T[Z\backslash V(\widetilde{P}_{h})]}(x_{\pi^{-1}(i)})$ and $N^{-}_{T[Z\backslash V(\widetilde{P}_{h})]}(y_{i})$, respectively.
Then
$$|X^{'}|\geq (n-3k+\pi^{-1}(i))/2-2h\geq (n-6.5k+5)/2\geq k-1,$$ and
$$|Y^{'}|\geq (n-3k+i)/2-2h\geq (n-6.5k+4)/2\geq k-1.$$
By the maximality of $h$, we also have $X^{'}\cap Y^{'}=\emptyset$ and $A(X^{'}\cup \{x_{\pi^{-1}(i)}\},Y^{'}\cup \{y_{i}\})=\emptyset$. Note that $|X^{'}\cup \{x_{\pi^{-1}(i)}\}|\geq k$ and $|Y^{'}\cup \{y_{i}\}|\geq k$. We can choose $Y^{''}\subseteq Y^{'}\cup \{y_{i}\}$ and $X^{''}\subseteq X^{'}\cup \{x_{\pi^{-1}(i)}\}$ such that $|X^{''}|=|Y^{''}|=k$ and $Y^{''}$ anchors $X^{''}$.

\vspace{0.2cm}

This completes the proof of the lemma.



\section{Proof of Theorem~\ref{main theorem}}

\noindent

Let $T$ be a $\lceil 12.5k-6\rceil$-connected tournament with $\delta^{+}(T)\geq 21k-14$. For $k=1$, the theorem is trivial. For $k=2$, the theorem follows directly from the fact that $g(2)=5$ due to Thomassen~\cite{T84}. Now assume $k\geq 3$.
 Let $X_{0}=\{x_{1},x_{2},\ldots,x_{k}\}$ and $Y_{0}=\{y_{1},y_{2},\ldots,y_{k}\}$ be two arbitrary disjoint vertex sets in $V(T)$. Our goal is to construct mutually disjoint dipaths $P_{x_{1}y_{1}},P_{x_{2}y_{2}},\ldots,P_{x_{k}y_{k}}$.
	
Let $T_{0}=T-(X_{0}\cup Y_{0})$. Choose $u_{1}\in V(T_{0})$ with $d^{+}_{T_{0}}(u_{1})=\delta^{+}(T_{0})$ and $v_{1}\in N^{+}_{T_{0}}(u_{1})$ with $d^{+}_{\small T[N^{+}_{T_{0}}(u_{1})]}(v_{1})=\delta^{+}(T[N^{+}_{T_{0}}(u_{1})])$. Indeed, $u_{1}$ and $v_{1}$ exist as $\delta^{+}(T_{0})\ge \delta^+(T)-2k\ge 19k-14>0$. Let $D_{1}=\{u_{1},v_{1}\}$ and $T_1=T_0-D_1$. Set $A_{1}=N^{+}_{T_{0}}(u_{1})\cap N^{+}_{T_{0}}(v_{1})$.
Let $k^{*}=\lceil 8.5k-6\rceil$.
For all $2\le i\le k^{*}$, repeating the above process, we can recursively construct $T_{i}=T_{i-1}-D_{i}$, where $D_{i}=\{u_{i},v_{i}\}$ such that $u_{i}\in V(T_{i-1})$ with $d^{+}_{T_{i-1}}(u_{i})=\delta^{+}(T_{i-1})$ and $v_{i}\in N^{+}_{T_{i-1}}(u_{i})$ with $d^{+}_{T[N^{+}_{T_{i-1}}(u_{i})]}(v_{i})=\delta^{+}(T[N^{+}_{T_{i-1}}(u_{i})])$.
Set $A_{i}=N^{+}_{T_{i-1}}(u_{i})\cap N^{+}_{T_{i-1}}(v_{i})$.
The vertices $u_{i}$ and $v_{i}$ exist since $\delta^{+}(T_{i-1})\geq 21k-14-2k-2k^*\ge 2k-3>0$.
\begin{claim}\label{CL: |A_i|}
For $i\in[k^*]$,	$|A_{i}|\leq\frac 12 (d^+_{T_{i-1}}(u_i)-1).$
\end{claim}
\begin{proof}[Proof of the claim:]
By the definition of $A_i$,  $$|A_{i}|=d^{+}_{T[N^{+}_{T_{i-1}}(u_{i})]}(v_{i}) = \delta^{+}(T[N^{+}_{T_{i-1}}(u_{i})]).$$
By Proposition~\ref{fact}, we have
\begin{align*}
	\delta^{+}(T[N^{+}_{T_{i-1}}(u_{i})])
	&\leq \frac 12\left(|N^{+}_{T_{i-1}}(u_{i})|-1\right)=\frac 12\left(d^+_{T_{i-1}}(u_i)-1\right).
\end{align*}
Therefore, we obtain that
\begin{align*}
	|A_{i}|\leq \frac 12\left(d^+_{T_{i-1}}(u_i)-1\right).
\end{align*}
\end{proof}	



Denote $T^{*} = T_{k^{*}}$, i.e.,  $T^{*}=T-(X_{0}\cup Y_{0}\cup D_{1}\cup D_{2}\cup \ldots \cup D_{k^{*}})$. Let $U = \{u_{1},u_{2},\ldots,u_{k^{*}}\}$ and $V = \{v_{1},v_{2},\ldots,v_{k^{*}}\}$. Note that $|V|=k^*=\lceil 8.5k-6\rceil$. Applying Lemma~\ref{key lemma} to the tournament $T[V]$, we have two sets $V_{1}= \{v_{\alpha_1},v_{\alpha_2},\ldots,v_{\alpha_k}\}$ and $V_2 = \{v_{\beta_1},v_{\beta_2},\ldots,v_{\beta_k}\}$ such that $V_{1}\cap V_{2}=\emptyset$, and $V_{1}$ anchors $V_{2}$ in $T[V]$.

For any $i\in [k]$, we have $N_{T^{*}}^{+}(x_{i})=N_{T}^{+}(x_{i})\backslash (X_0\cup Y_0\cup U\cup V)$. It follows that
$$d_{T^{*}}^{+}(x_{i})\geq d_{T}^{+}(x_{i})-|X_{0}\cup Y_{0}\cup U\cup V|\geq (21k-14)-(2k^{*}+2k)\geq 2k-3\geq k.$$
Hence, we can greedily find $X_{1}=\{x_{1}',x_{2}',\ldots,x_{k}'\} \subseteq V(T^{*})$ such that $x_{i}\rightarrow x_{i}'$ and for all $i\in[k]$.

We shrink every $D_{\alpha_{i}}=\{u_{\alpha_{i}},v_{\alpha_{i}}\}$ to a single vertex $z_{\alpha_{i}}$  and let $Z=\{z_{\alpha_{1}},z_{\alpha_{2}},\ldots,z_{\alpha_{k}}\}$. Now we construct an auxiliary bipartite graph $G=(X_{1}\cup Z, E)$ such that $(x_{i}',z_{\alpha_{j}})\in E(G)$ if and only if $x_{i}'$ dominates at least one of $\{u_{\alpha_{j}},v_{\alpha_{j}}\}$. Let $M$ be a maximum matching of $G$. By Lemma~\ref{matching}, there is a subset $S\subseteq X_{1}$ such that $|M|=|X_{1}|-(|S|-|N_{G}(S)|)$. Let $d=|S|-|N_{G}(S)|$. Then $0\leq d\leq |S|=s$. We may assume, without loss of generality, that $S=\{x_{1}',x_{2}',\ldots,x_{s}'\}$, $N_{G}(S)=\{z_{\alpha _{d+1}},z_{\alpha_{d+2}},\ldots,z_{\alpha_{s}}\}$ and $M=\{x_{d+1}'z_{\alpha_{d+1}},x_{d+2}'z_{\alpha_{d+2}},\ldots,x_{k}'z_{\alpha_{k}}\}$.
Then, for $i\in[d+1,k]$, by the definition of $G$, at least one of $x_{i}^{'}\rightarrow u_{\alpha_{i}}\rightarrow v_{\alpha_{i}}$ and $x_{i}^{'}\rightarrow v_{\alpha_{i}}$ must exist. This indicates that we can find $k-d$ disjoint dipaths $Q_{d+1}, \dots, Q_k$ such that $Q_i$ is from $x_{i}$ to $v_{\alpha_{i}}$ for $i\in [d+1, k]$.
If $d=0$, then we have found $k$ disjoint dipaths $Q_i$ from $x_{i}$ to $v_{\alpha_{i}}$ for every $i\in [k]$.

Now assume $d\geq 1$. Note that, for each $x_i'\in S$, $x_i'z_{\alpha_{j}}\notin E(G)$ for $j\in[d]\cup [s+1,k]$. By the construction of $G$,  $A_T(x_{i}', \{u_{\alpha_{j}},v_{\alpha_{j}}\})=\emptyset$ for every $j\in [d]\cup [s+1, k]$. Since $T$ is a tournament and $x_i'\in V(T^{*})\subseteq V(T_{\alpha_j-1})$, we have $x_i'\in N^+_{T_{\alpha_j-1}}(u_{\alpha_j})\cap N^+_{T_{\alpha_j-1}}(v_{\alpha_j})=A_{\alpha_j}$ for   $j\in [d] \cup [s+1, k]$. This follows that $S\subseteq A_{\alpha_{j}}$ for all $j\in[d]\cup [s+1,k]$.

In the rest of our proof, we will construct $d$ disjoint dipaths $Q_{1}, \dots, Q_d$ such that $Q_i$ is from $x_{i}$ to $v_{\alpha_{i}}$ for $i\in [d]$.

For each $i\in [d]$, let 
$$F_{i}=X_{0}\cup Y_{0}\cup \left(\cup_{j=1}^{\alpha_{i}-1}D_{j}\right)\cup \widetilde{U}\cup \widetilde{V}\cup X_{1}\cup A_{\alpha_{i}},$$ where $\widetilde{U}=\{u_{\alpha_{d+1}},u_{\alpha_{d+2}},\ldots,u_{\alpha_{s}}\}$ and $\widetilde{V}=V\backslash \{v_{1},v_{2},\ldots,v_{\alpha_{i}}\}$.
We will construct a dipath $Q_i$ that does not pass by the vertices of $F_i$. 
We first estimate $d^{+}_{T}(x_{i}', F_{i})$. It is clear that
 \begin{align}
	d^{+}_{T}(x_{i}', F_{i})
&\leq d^{+}_{T}(x_{i}', X_{0}\cup Y_{0})+d^{+}_{T}(x_{i}', \cup_{j=1}^{\alpha_{i}-1}D_{j})+d^{+}_{T}(x_{i}',\widetilde{U}\cup \widetilde{V})+d^{+}_{T}(x_{i}', X_{1}\cup A_{\alpha_{i}}).\label{EQ:e1}
\end{align}
For convenience, let $d^{+}_{T}(x_{i}', X_{0}\cup Y_{0})=\gamma_{i}$ and $d^{+}_{T}(x_{i}', \cup_{j=1}^{\alpha_{i}-1}D_{j})=\tau_{i}$, respectively. As $x_{i}\rightarrow x_{i}'$, we have
\begin{align}
	d^{+}_{T}(x_{i}', X_{0}\cup Y_{0})=\gamma_{i}\leq 2k-1.\label{EQ:e2}
\end{align}
Clearly,
\begin{align}
	d^{+}_{T}\left(x_{i}', \cup_{j=1}^{\alpha_{i}-1}D_{j}\right)=\tau_{i}\leq \sum_{j=1}^{\alpha_i-1}|D_j|=2(\alpha_{i}-1),\label{EQ:e3}
\end{align}
and 
\begin{align}
	d^{+}_{T}(x_{i}',\widetilde{U}\cup \widetilde{V})\le |\widetilde{U}|+|\widetilde{V}|= s-d+k^{*}-\alpha_{i}.\label{EQ:e4}
\end{align}
Recall that $S\subseteq A_{\alpha_{i}}$ for any $i\in [d]$. Therefore, 
\begin{align}
	d^{+}_{T}(x_{i}', X_{1}\cup A_{\alpha_{i}})
	&=d^{+}_{T}(x_{i}', (X_{1}\backslash S)\cup A_{\alpha_{i}})\nonumber\\
	&\leq d^{+}_{T}(x_{i}', (X_{1}\backslash S))+ d^{+}_{T}(x_{i}', A_{\alpha_{i}})\nonumber\\
	&\leq |X_{1}\backslash S|+|A_{\alpha_{i}}|\nonumber\\
	&\leq k-s+\frac 12\left(d^{+}_{T_{\alpha_{i}-1}}(u_{\alpha_{i}})-1\right)\nonumber\\
	&\leq k-s+\frac 12\left(d^{+}_{T}(x_{i}')-\gamma_{i}-\tau_{i}-1\right),\label{EQ:e5}
\end{align}
where the third inequality holds by Claim~\ref{CL: |A_i|}, and the fourth inequality holds since \begin{align*}
d^{+}_{T_{\alpha_{i}-1}}(u_{\alpha_{i}})&=\delta^{+}(T_{\alpha_{i}-1})\leq d^{+}_{T_{\alpha_{i}-1}}(x_{i}')\\
&=d^{+}_{T}(x_{i}')-d^{+}_{T}(x_{i}', X_{0}\cup Y_{0})-d^{+}_{T}(x_{i}', \cup_{j=1}^{\alpha_{i}-1}D_{j})\\
&=d^{+}_{T}(x_{i}')-\gamma_{i}-\tau_{i}.
\end{align*}
Substitute (\ref{EQ:e2}),(\ref{EQ:e3}),(\ref{EQ:e4}) and (\ref{EQ:e5}) into (\ref{EQ:e1}), for every $i\in [d]$, we have
 \begin{align*}
	d^{+}_{T}(x_{i}', F_{i})
&\leq \gamma_{i}+\tau_{i}+s-d+k^{*}-\alpha_{i}+k-s+\frac 12\left(d^{+}_{T}(x_{i}')-\gamma_{i}-\tau_{i}-1\right)\\
&=\frac 12 d^{+}_{T}(x_{i}')+\frac 12\left(\gamma_{i}+\tau_{i}-1\right)+k^{*}+k-\alpha_{i}-d\\
&\leq \frac 12 d^{+}_{T}(x_{i}')+\frac 12\left(2k-1+2(\alpha_{i}-1)-1\right)+k^{*}+k-\alpha_{i}-d\\
&= \frac 12 d^{+}_{T}(x_{i}')+k^{*}+2k-d-2,
\end{align*}
where the second inequality holds since $\gamma_{i}\leq 2k-1$ and $\tau_{i}\leq 2(\alpha_{i}-1)$.
Therefore,  
 \begin{align*}
 d^{+}_{T}(x_{i}', V(T)\backslash F_{i})
 &=d^{+}_{T}(x_{i}^{'})-d^{+}_{T}(x_{i}', F_{i})\\
&\geq \frac 12 d^{+}_{T}(x_{i}')-k^{*}-2k+d+2\\
&\geq \frac 12 \left(21k-14\right)-(8.5k-5.5)-2k+d+2\\
&>d.
\end{align*}
where the second inequality holds since $k^{*}=\lceil 8.5k-6\rceil\leq 8.5k-5.5$.

Hence, we can greedily find $d$ vertices $x_{1}'', x_{2}'', \ldots, x_{d}''$ such that, for every $i\in [d]$, $x_{i}'\rightarrow x_{i}''$ and $x_{i}''\in V(T)\backslash F_{i}$. Then, $x_{i}''$ dominates at least one of $u_{\alpha_{i}}$ and $v_{\alpha_{i}}$. This indicates that at least one of $x_{i}\rightarrow x_{i}^{'}\rightarrow x_{i}^{''}\rightarrow u_{\alpha_{i}}\rightarrow v_{\alpha_{i}}$, $x_{i}\rightarrow x_{i}^{'}\rightarrow x_{i}^{''}\rightarrow v_{\alpha_{i}}$ must exist.
Therefore, we can find $d$ disjoint dipaths $Q_1, \dots, Q_d$ such that $Q_i$ is from $x_i$ to $v_{\alpha_i}$ for $i\in[d]$.

Note that $$\bigcup_{i=1}^dV(Q_i)\subseteq\left(\bigcup_{i=1}^d\{x_i, x_i',x_i''\}\right)\cup\left(\bigcup_{i=1}^d D_{\alpha_i}\right)$$ and recall that $$\bigcup_{i=d+1}^kV(Q_i)\subseteq\left(\bigcup_{i=d+1}^k\{x_i, x_i'\}\right)\cup\left(\bigcup_{i=d+1}^k D_{\alpha_i}\right).$$
We obtain that $Q_{1},\dots,Q_d, Q_{d+1},\ldots, Q_{k}$ are $k$ disjoint dipaths  such that $Q_{i}$ is from $x_i$ to $v_{\alpha_i}$ and every one has length at most $4$ for $i\in [k]$.

Let $B=\left(\bigcup_{i=1}^kV(Q_{i})\right)\cup \left(V\backslash (V_{1}\cup V_{2})\right)$. Note that $|V(Q_{i})|\leq 5$. Thus $$|B|\leq 5k+( k^*-2k)\leq \lceil 11.5k-6\rceil.$$
Since $T$ is $\lceil 12.5k-6\rceil$-connected, we have $T-B$ is $k$-connected. By Corollary~\ref{corollary}, there are $k$ disjoint dipaths $R_{i}^{\pi}$ from $v_{\beta_{i}}$ to $y_{\pi(i)}$ in $T-B$ for some permutation $\pi$ of $[k]$. Due to $V_{1}$ anchors $V_{2}$ in $T[V]$, there exist $k$ disjoint dipaths $P_{i}^{\pi^{-1}}$ from $v_{\alpha_{i}}$ to $v_{\beta_{\pi^{-1}(i)}}$ in $T[V]$.

Finally, let $P_{x_{i}y_{i}}=Q_i\cup P_i^{\pi^{-1}}\cup R_{i}^{\pi}$ for $i\in[k]$. Then $P_{x_{1}y_{1}},P_{x_{2}y_{2}},\ldots,P_{x_{k}y_{k}}$ are the desired dipaths.

This completes our proof of the theorem.

\section{Remarks and Discussions}

In this paper, we prove that every $\lceil 12.5k-6\rceil$-connected tournament with minimum out-degree at least $21k-14$ is $k$-linked.
{In fact, authors in~\cite{BJJ22} and~\cite{ZhouQiYan23} attempted to give  `better' bounds of the minimum out-degree. However, there is a similar  bug in their proofs. Specifically, in page 3 of~\cite{BJJ22}, since vertices in $U$ are relabeled when choosing $U'$, for some $u_i'$ in $U'$, it could happen that some $x_i''=u_j\in U^*$ is not in $A_i'$ (nor $X_0, Y_0, U', V, X'$), but an out-neighbor of both $u_i'$ and $v_i'$.}


  It will be very interesting to determine optimal values of $f(k)$ and $d(k)$ to guarantee a digraph  $k$-linked  under both  connectivity $f(k)$ and minimum out-degree $d(k)$  are linear $k$.
We believe that $8.5k-6$ is not optimal in the anchor lemma (Lemma~\ref{key lemma}), what is the exact value there.

\subsection*{Acknowledgement}

\noindent

The work of the first, second and fourth author was supported by the National Natural Science Foundation of China (No. 12071453), the National Key R and D Program of China(2020YFA0713100),  and the Innovation Program for Quantum Science and Technology, China (2021ZD0302902).

\vskip 3mm
\end{spacing}
\end{document}